\documentclass[12pt, a4paper]{article}
\usepackage[pagewise]{lineno}
\usepackage{graphicx}
\usepackage{amssymb}
\usepackage{latexsym, bm}
\usepackage{multicol}
\usepackage{indentfirst}
\usepackage{amssymb,amsfonts}
\usepackage{amsmath}
\usepackage{cases}
\usepackage[noadjust]{cite}
\usepackage[colorlinks,citecolor=red]{hyperref}
\newtheorem{theorem}{Theorem}
\newtheorem{lemma}{Lemma}

\newtheorem{definition}{Definition}

\newtheorem{remark}{Remark}

\usepackage{amssymb}
\numberwithin{equation}{section}
\numberwithin{theorem}{section}
\numberwithin{remark}{section}
\numberwithin{definition}{section}
\numberwithin{lemma}{section}
\numberwithin{corollary}{section}
\numberwithin{proposition}{section}
\numberwithin{notation}{section}
\title{New Liouville type theorems for the stationary MHD
equations in $\mathbb{R}^3$}
\author{Wenke Tan\footnote{tanwenkeybfq@163.com}\\
{\small Key Laboratory of Computing and Stochastic Mathematics (Ministry of Education),}\\
{\small School of Mathematics and Statistics, Hunan Normal University,}\\
{\small Changsha, Hunan 410081, China}\\
}

\date{}

\begin{document}
\maketitle
{\bf Abstract:} We research the Liouville type problem
for the 3D stationary MHD equations in the frequency space. We establish two new Liouville type theorems for solutions with finite Dirichlet energy. Specifically, we show that the low-frequency part of the velocity field plays the leading role in a Liouville theory for MHD equations and then improve the results of Chae-Weng \cite{Chae-W}.

\medskip

\section{Introduction}
In this paper, we consider the homogeneous stationary Magnetohydrodynamics (MHD)
equations in the whole space $\mathbb{R}^3$:
\begin{equation}\label{MHD}
 \left\{\begin{array}{ll}
-\Delta u+u\cdot\nabla u+\nabla P=B\cdot\nabla B,\\
-\Delta B+u\cdot\nabla B=B\cdot\nabla u,\\
\nabla\cdot u=\nabla\cdot B=0,
\end{array}\right.
\end{equation}
where the unknowns $u(x)$ and $B(x)$ stand
for the velocity vector and the magnetic field respectively, and $P(x)$ denotes the scalar
pressure. This system and its time-dependent analogue are used to model electrically
conductive fluids such as plasma, liquid metals, electrolytes etc. 
For more physical background we refer to Schnack \cite{Schnack} and the
references therein.

When $B=0$, the equations are reduced to the Navier-Stokes equations. It is well known
that due to the pioneering work of Leray \cite{Leray}, it has been an open problem whether $u=0$ is the only solution in the class of the finite Dirichlet integral $$D(u)=\int_{\mathbb{R}^3}|\nabla u|^2dx$$ with the homogeneous condition at infinity
$$\lim_{|x|\to\infty}u(x)=0.$$
This is a famous Liouville type statement on the stationary Navier-Stokes equations. For the 2D Navier-Stokes
equations, Gilbarg and Weinberger proved the Liouville-type theorem in \cite{G-W}. The
authors used the idea that the vorticity $\omega=\nabla\times u$ is a scalar satisfying an elliptic equation that
enables one to apply the maximum principle. This, together with another result showing
that $\omega\to0$ at infinity, implies that $\omega=0$. From this, the authors deduced that $u$ and
$P$ are constant in $\mathbb{R}^2$. The same approach as in \cite{G-W} fails in the 3D case, mainly due to
the more complicated form of the equation for vorticity.  Regarding the 3D case, one of
the best attempts made to solve the above or related problems was presented by Galdi
in \cite{Galdi}, where he showed that if $u\in L^\frac{9}{2}(\mathbb{R}^3)$, then it holds that $u=0$.  In the case of the
n-dimensional Navier-Stokes system with $n\geq 4$, the problem was resolved by Galdi in \cite{Galdi}
by similar argument to the case of $u\in L^\frac{9}{2}(\mathbb{R}^3)$ with additional assumption $\lim_{|x|\to\infty}|u(x)|=0$. It is worth pointing out that when $n\geq 5$, the Liouville
problem with only finite Dirichlet integral is still open. Chae and Wolf \cite{C-W} gave a logarithmic improvement of Galdi's
result by assuming that
\begin{align*}
\int_{\mathbb{R}^3}|u|^\frac{9}{2}\{\ln(2+|u|^{-1})\}^{-1}dx<\infty.
\end{align*} Also, Chae \cite{C} showed that the condition
\begin{align*}
\Delta u\in L^\frac{6}{5}(\mathbb{R}^3)
\end{align*} is sufficient for $u=0$ in $\mathbb{R}^3$. Kozono, Terasawa and Wakasugi proved in \cite{Kozono} that $u=0$ if the vorticity satisfies
\begin{align*}
\lim_{|x|\to \infty}|x|^\frac{5}{3}|\omega(x)|\leq(\delta D(u))^\frac{1}{3}
\end{align*} or the velocity satisfies
\begin{align*}
||u||_{L^{\frac{9}{2},\infty}(\mathbb{R}^3)}\leq (\delta D(u))^\frac{1}{3}
\end{align*} for a small constant $\delta$. Then, the restriction imposed on the norm $||u||_{L^{\frac{9}{2},\infty}(\mathbb{R}^3)}$ in \cite{Kozono} was relaxed by Seregin and Wang in \cite{S-W}. For more references on Liouville type theorem
for stationary Navier-Stokes equations, we refer to \cite{C-J-L,C-N-Y,J,S2,S,Tsai} and the references
therein.

If we take the magnetic effect into account, it is natural to consider the corresponding
Liouville type problem for solutions $(u,B)$ of \eqref{MHD} in the class of the finite Dirichlet
integral
\begin{align}\label{DI}
D(u,B)=\int_{\mathbb{R}^3}(|\nabla u|^2+|\nabla B|^2)dx<\infty
\end{align}with the condition
\begin{align}\label{far}
\lim_{|x|\to\infty}|u(x)|=\lim_{|x|\to\infty}|B(x)|=0.
\end{align}However, for the MHD system, the situation is quite different. The involvement of the
magnetic field makes the problem much more complicated. For the 2D MHD equations,
Wang-Wang \cite{W-W} have established the Liouville type theorems by assuming the smallness of
the norm of the magnetic field. The smallness was removed by Nicola-Francis-Simon \cite{N-F-S}.
For the 3D MHD equations, in \cite{Sch}, Schulz proved that if the smooth solution $(u,B)$ of
the stationary MHD equations \eqref{MHD} is in $L^6(\mathbb{R}^3)\cap BMO^{-1}(\mathbb{R}^3)$ then it is identically zero.
Chae-Wolf \cite{C-W2} showed that $L^6$ mean oscillations of the potential function of the velocity
and magnetic field have certain linear growth by using the technique of \cite{C-W}. Recently,
Li-Pan \cite{L-P} proved two forms of Liouville theorems for D-solutions, one is the case where
$u\to(1,0,0)$ and $B\to(0,0,0)$ for any viscosity and magnetic resistivity; another case is
$u\to(0,0,0)$ and $B\to(1,0,0)$ by taking the same viscosity and magnetic resistivity. This
result was improved by Wang-Yang \cite{W-Y}, they showed that $u$ and $B$ are constant for
the MHD system when the velocity field tends to a constant vector at infinity while the
magnetic field tends to $0$ without any assumptions on viscosity, magnetic resistivity. On
the other hand, some numerical experiments in \cite{P-P-S} seem to indicate that the velocity field
should play a more important role than the magnetic field in the regularity theory. See also
\cite{H-X,H-X2,W-Z}. One may wonder whether the velocity fields also play the leading role in a Liouville
theory for MHD equations. Partial progress has been made for the three-dimensional case,
where Chae-Weng \cite{Chae-W} recently proved the axially symmetric solution $(u,B)$ with finite
Dirichlet integral is trivial if $u\in L^3(\mathbb{R}^3)$.  For more references on Liouville type theorem
for stationary MHD equations, we refer to \cite{F-W,L-L-N,Y-X} and the references therein.

To the author's knowledge, those authors of the previous works in the literature are analyzing the Liouville problem in physical space. There are no related results that consider
the Liouville problem in frequency space. In this paper, we will exploit the Littlewood-Paley theory to research the Liouville problem for the stationary MHD equations \eqref{MHD}
with assumptions \eqref{DI}. Our analysis starts from the scaling symmetry $(u(x),B(x),P(x))\rightarrow(u_\lambda(x),B_\lambda(x),P_\lambda(x))$ of stationary MHD equations \eqref{MHD} for $\lambda>0$, where
 \begin{align*}
 u_\lambda(x)\doteq \lambda u(\lambda x)\quad B_\lambda(x)\doteq \lambda B(\lambda x)\quad P_\lambda(x)\doteq \lambda^2 P(\lambda x).
 \end{align*}
 Among other things, this means the quantity $\int_{\mathbb{R}^3}(|\nabla u|^2+|\nabla B|^2)dx$ is subcritical. This character inspires us that the main obstacle to solving the Liouville problem may come from the low-frequency parts of solutions. To exploit the subcritical character of \eqref{DI}, we first use the localization techniques in frequency space to transform the Dirichlet integral to a new integral which is completely determined by the behavior of the solution $(u,B)$ in any neighborhood of the origin in frequency space. Based on this transformation, we establish two new Liouville-type theorems.

Before proceeding with our main results, we define the weak solution to the MHD system \eqref{MHD}.
\begin{definition}
We say that $(u,B)\in\mathcal{S}'(\mathbb{R}^3)\times\mathcal{S}'(\mathbb{R}^3)$ is a weak solution to the MHD
equations \eqref{MHD} in $\mathbb{R}^3$ provided that:\\
(1) $\nabla u, \nabla B\in L^2_{loc}(\mathbb{R}^3)$;\\
(2) $\nabla\cdot u=\nabla\cdot B=0$ in the sense of distributions;\\
(3) The following identity holds
\begin{align*}
\int_{\mathbb{R}^3}\nabla u:\nabla \varphi dx=&-\int_{\mathbb{R}^3}u\cdot\nabla u\cdot \varphi dx+\int_{\mathbb{R}^3}B\cdot\nabla B\cdot \varphi dx,\\
\int_{\mathbb{R}^3}\nabla B:\nabla \varphi dx=&-\int_{\mathbb{R}^3}u\cdot\nabla B\cdot \varphi dx+\int_{\mathbb{R}^3}B\cdot\nabla u\cdot \varphi dx
\end{align*}for all $\varphi(x)\in (C_0^\infty(\mathbb{R}^3))^3$ with $\nabla\cdot\varphi=0$.
\end{definition}

Our first main result is as follows,
\begin{theorem}\label{main}
Let $(u,B)$ be a weak solution of \eqref{MHD} in the class \eqref{DI}. Then
$u=B=0$ if the following condition is valid
\begin{align}\label{cond1}
\liminf_{k\to-\infty}2^{-k}(||\dot{S}_ku||_{L^\infty(\mathbb{R}^3)}+||\dot{S}_kB||_{L^\infty(\mathbb{R}^3)})<\infty.
\end{align}Specifically, we deduce that if $(u,B)$ satisfies
\begin{align}\label{cond2}
\liminf_{k\to-\infty}(||\dot{S}_ku||_{\dot{B}^{-1}_{\infty,\infty}}+||\dot{S}_kB||_{\dot{B}^{-1}_{\infty,\infty}})<\infty,
\end{align}then $u=B=0$.
\end{theorem}
\begin{remark}
Noticing that
\begin{align*}
supp\mathcal{F}(\dot{S}_kf)(\xi)\subset\{\xi:|\xi|<2^k\}\to \{0\}~as~k\to-\infty,
\end{align*}Theorem \ref{main} implies that the uniqueness of the solution $(u,B)$ is completely determined by the local information of $(u,B)$ at the origin in frequency space.
\end{remark}

Our second result shows that the low-frequency part of velocity fields plays the leading role in a Liouville theory for MHD equations \eqref{MHD}.
\begin{theorem}\label{main1}
Let $(u,B)$ be a weak solution of \eqref{MHD} in the class \eqref{DI}. If the following condition is valid
\begin{align}\label{cond3}
\liminf_{k\to-\infty}||\dot{S}_{k}u||_{L^3(\mathbb{R}^3)}<\infty,
\end{align}
then
$u=B=0$.
\end{theorem}
\begin{remark}
Chae-Weng \cite{Chae-W} recently proved the axially symmetric solution $(u,B)$ with finite
Dirichlet integral is trivial if $u\in L^3(\mathbb{R}^3)$. We remove the restriction of axial symmetry and show that the uniqueness of solutions is completely determined by the $L^3$ norm of the low-frequency part near the origin.
\end{remark}

Our paper is organized as follows: In Section 2, we list some notations and recall some
Lemmas about the Littlewood-Paley theory which will be used in the sequel. In Section 3 we present the proofs of results.

\section{Preliminaries}

{\bf Notation:} In this paper, we denote $B_r=\{x\in\mathbb{R}^3:|x|\leq r\}$.  The matrix $\nabla u$ denotes the gradient of $u$ with respect to the $x$ variable, whose
$(i,j)$-th component is given by $(\nabla u)_{i,j}=\partial_ju_i$ with $1\leq i,j\leq3$. Throughout this
paper, $C$ stands for some positive constant, which may vary from line to line. Given a Banach space $X$, we denote its norm by $||\cdot||_{X}$.
We use $\mathcal{S}(\mathbb{R})^3$ and $\mathcal{S}'(\mathbb{R}^3)$ to denote Schwartz functions and the tempered distributions spaces on $\mathbb{R}^3$, respectively.

Next, we will recall some well-known facts about the Littlewood-Paley decomposition. Firstly, let us recall that for every $f\in\mathcal{S}'(\mathbb{R}^3)$
the Fourier transform of $f$ is defined by
\begin{align*}
(\mathcal{F}f)(\xi)=\hat{f}(\xi)=(2\pi)^{-\frac{3}{2}}\int_{\mathbb{R}^3}e^{-ix\cdot\xi}f(x)dx\quad for\quad\xi\in\mathbb{R}^3
\end{align*}The inverse Fourier transform of $g\in\mathcal{S}'(\mathbb{R}^3)$ is given by
\begin{align*}
(\mathcal{F}^{-1}g)(x)=\breve{g}(x)=(2\pi)^{-\frac{3}{2}}\int_{\mathbb{R}^3}e^{ix\cdot\xi}g(\xi)dx\quad for \quad x\in\mathbb{R}^3.
\end{align*}By using Fourier transform, we can define homogeneous Sobolev norm $||\cdot||_{\dot{H}^s}$ with $s\in\mathbb{R}$ as
\begin{align*}
||f||_{\dot{H}^s(\mathbb{R}^3)}=(\int_{\mathbb{R}^3}|\xi|^{2s}|\hat{f}(\xi)|^2d\xi)^\frac{1}{2}.
\end{align*}

The main tool in this paper is the Littlewood-Paley decomposition of distributions into dyadic blocks of frequencies:
\begin{definition}\label{LPD}
Let $\psi(\xi)\in C^\infty_0(B_1)$ be a non-negative function so that $\psi(\xi)=1$ for $|\xi|\leq \frac{1}{2}$. Let $\varphi(\xi)$ be defined as $\varphi(\xi)=\psi(2^{-1}\xi)-\psi(\xi)$.
For given $u\in\mathcal{S}'(\mathbb{R}^3)$, the homogeneous dyadic blocks $\dot{\Delta}_k$ and the homogeneous low-frequency cut-off operator $\dot{S}_k$ are defined for all $k\in\mathbb{Z}$ by
\begin{align}\label{block}
\dot{\Delta}_ku(x)&=\mathcal{F}^{-1}(\varphi(2^{-k}\cdot)\hat{u}(\cdot))(x)=\frac{1}{(2\pi)^\frac{3}{2}}\int_{\mathbb{R}^3}e^{ix\cdot\xi}\varphi(2^{-k}\xi)\hat{u}(\xi)d\xi,\\
\dot{S}_ku(x)&=\mathcal{F}^{-1}(\psi(2^{-k}\cdot)\hat{u}(\cdot))(x)=\frac{1}{(2\pi)^\frac{3}{2}}\int_{\mathbb{R}^3}e^{ix\cdot\xi}\psi(2^{-k}\xi)\hat{u}(\xi)d\xi.
\end{align}
\end{definition}

 Throughout this paper, we will use the notation that $\tilde{\dot{\Delta}}_lu=\sum_{|l'-l|\leq2}\dot{\Delta}_{l'}u$. In the homogeneous case, the following Littlewood-Paley decomposition makes sense
\begin{align*}
u(x)=\sum_{k\in\mathbb{Z}}\dot{\Delta}_ku(x)\quad for\quad u\in\mathcal{S}_h'(\mathbb{R}^3),
\end{align*} where $\mathcal{S}_h'(\mathbb{R}^3)$ is given by
\begin{align*}
\mathcal{S}_h'(\mathbb{R}^3)=\{u\in\mathcal{S}'(\mathbb{R}^3):\lim_{k\to-\infty}\dot{S}_ku=0\}.
\end{align*}Moreover, it holds that
\begin{align*}
\dot{S}_ku(x)=\sum_{l\leq k-1}\dot{\Delta}_lu(x).
\end{align*} Based on the Littlewood-Paley decomposition, we introduce the homogeneous Besov space $\dot{B}^s_{p,q}$ as follows
\begin{definition}
Let $s\in\mathbb{R}$ and $1\leq p,q\leq\infty$. The homogeneous Besov space $\dot{B}^s_{p,q}$ consists of all tempered distribution $f$ such that
\begin{align*}
\dot{B}^s_{p,q}=\{f\in\mathcal{S}'_h:\|f\|_{\dot{B}^s_{p,q}}<\infty\},
\end{align*}
where
\begin{align*}
		\begin{split}
			\| f\|_{\dot{B}^s_{p,q}}= \left \{
			\begin{array}{ll}
				\left(\sum\limits_{j\in\mathbb{Z}}2^{sjq}\|\dot{\Delta}_jf\|_{L^p}^q\right)^{\frac{1}{q}}                   & if \ 1\leq q <\infty,\\
				\sup\limits_{j\in\mathbb{Z}}2^{sj}\|\dot{\Delta}_jf\|_{L^p}                                & if \ q=\infty.
			\end{array}
			\right.
		\end{split}
	\end{align*}

\end{definition}

A great advantage of the localized techniques in frequency is the so-called Bernstein inequalities which will be used in the sequel.
\begin{lemma}\cite{B-C-D}
Let $\mathcal{C}$ be an annulus and $B$ a ball. A constant C exists such that
for any nonnegative integer $k$, any couple $(p,q)$ in $[1,\infty]^2$ with $1\leq p\leq q$, and any function $u\in L^p$, we have
\begin{align}
supp\hat{u}\subset\lambda B&\Rightarrow ||\nabla^k u||_{L^q}\leq C^{k+1}\lambda^{k+\frac{3}{p}-\frac{3}{q}}||u||_{L^p},\label{Ber1}\\
supp\hat{u}\subset\lambda \mathcal{C}&\Rightarrow C^{-k-1}\lambda^k||u||_{L^p}\leq||\nabla^ku||_{L^p}\leq C^{k+1}\lambda^k||u||_{L^p}\label{Ber2}
\end{align}
\end{lemma}
From the Littlewood-Paley decomposition and Bernstein inequality \eqref{Ber2}, the homogeneous Sobolev norm $||\cdot||_{\dot{H}^s}$ can be equivalently written as
\begin{align*}
||f||_{\dot{H}^s(\mathbb{R}^3)}=(\sum_{k\in\mathbb{Z}}2^{2ks}||\dot{\Delta}_k f||^2_{L^2})^\frac{1}{2}.
\end{align*}

It is also worth pointing out that throughout this paper, we will use the following well-known facts:
\begin{align}\label{+10}
||\dot{S}_kf||_{L^p}\leq C||f||_{L^p},\quad ||f^k||_{L^p}\leq C||f||_{L^p}
\end{align}here $f^k=f-\dot{S}_kf$.

\section{Proofs of results}
Based on the Littlewood-Paley theory, we first derive a new formula for Dirichlet integral \eqref{DI}. The main novelty of this formula is that although the Dirichlet integral \eqref{DI}
is a globally defined quantity, but it can be characterized by the local information of $u$ and $B$ in any neighborhood of the origin in frequency space.
\begin{lemma}\label{lemma}
Let $(u,B)$ be a weak solution of \eqref{MHD} in the class \eqref{DI}, it follows
\begin{align}\label{formula}
&\int_{\mathbb{R}^3}(|\nabla u|^2+|\nabla B|^2)dx\\=&\liminf_{k\to-\infty}\{-\int_{\mathbb{R}^3}
\sum_{l=k-1}^{k+1}\sum_{l'=l-2}^{k-1}\dot{S}_{l-2}\dot{S}_{k}u\cdot\nabla\dot{\Delta}_{l'}u\cdot\dot{\Delta}_lu^kdx
-\int_{\mathbb{R}^3}\sum_{l=k-3}^{k}\dot{\Delta}_l\dot{S}_ku\cdot\nabla\dot{S}_ku\cdot\tilde{\dot{\Delta}}_lu^kdx\notag\\
&-\int_{\mathbb{R}^3}\sum_{l=k-1}^{k+1}\sum_{l'=l-2}^{k-1}\dot{S}_{l-2}\dot{S}_{k}u\cdot\nabla\dot{\Delta}_{l'}B\cdot\dot{\Delta}_lB^kdx
-\int_{\mathbb{R}^3}\sum_{l=k-3}^{k}\dot{\Delta}_l\dot{S}_ku\cdot\nabla\dot{S}_kB\cdot\tilde{\dot{\Delta}}_lB^kdx\notag\\
+&\int_{\mathbb{R}^3}\sum_{l=k-1}^{k+1}\sum_{l'=l-2}^{k-1}\dot{S}_{l-2}\dot{S}_{k}B\cdot\nabla\dot{\Delta}_{l'}B\cdot\dot{\Delta}_lu^kdx
+\int_{\mathbb{R}^3}\sum_{l=k-3}^{k}\dot{\Delta}_l\dot{S}_kB\cdot\nabla\dot{S}_kB\cdot\tilde{\dot{\Delta}}_lu^kdx\notag\\
&+\int_{\mathbb{R}^3}\sum_{l=k-1}^{k+1}\sum_{l'=l-2}^{k-1}\dot{S}_{l-2}\dot{S}_{k}B\cdot\nabla\dot{\Delta}_{l'}u\cdot\dot{\Delta}_lB^kdx
+\int_{\mathbb{R}^3}\sum_{l=k-3}^{k}\dot{\Delta}_l\dot{S}_kB\cdot\nabla\dot{S}_ku\cdot\tilde{\dot{\Delta}}_lB^kdx\notag\\
&-\int_{\mathbb{R}^3}\sum_{l=k-1}^{[\frac{3}{4} k]}\tilde{\dot{\Delta}}_{l}u^k\cdot\nabla\dot{S}_ku\cdot\dot{\Delta}_lu^kdx-\int_{\mathbb{R}^3}\sum_{l=k-1}^{[\frac{3}{4} k]}\tilde{\dot{\Delta}}_{l}u^k\cdot\nabla\dot{S}_kB\cdot\dot{\Delta}_lB^kdx\notag\\
&+\int_{\mathbb{R}^3}\sum_{l=k-1}^{[\frac{3}{4} k]}\tilde{\dot{\Delta}}_{l}B^k\cdot\nabla\dot{S}_kB\cdot\dot{\Delta}_lu^kdx+\int_{\mathbb{R}^3}\sum_{l=k-1}^{[\frac{3}{4} k]}\tilde{\dot{\Delta}}_{l}B^k\cdot\nabla\dot{S}_ku\cdot\dot{\Delta}_lB^kdx\notag\}.
\end{align}
\end{lemma}
\emph{\bf Proof of Lemma \ref{lemma}.}
We define $u^k(x)$ and $B^k(x)$ as follows
\begin{align}\label{hp}
u^k(x)\doteq\sum_{l\geq k}\dot{\Delta}_lu(x)=&u(x)-\dot{S}_ku(x),\\
B^k(x)\doteq\sum_{l\geq k}\dot{\Delta}_lB(x)=&B(x)-\dot{S}_kB(x).
\end{align}It is not difficult to see that $u^k$ and $B^k$ are obeying the following estimates
\begin{align*}
||\nabla u^k||_{L^2}&\leq ||\nabla u||_{L^2},\\
||u^k||_{L^3}\leq\sum_{l\geq k}||\dot{\Delta}_lu||_{L^3}\leq\sum_{l\geq k} 2^{\frac{l}{2}}||\dot{\Delta}_lu||_{L^2}&\leq\sum_{l\geq k}2^{\frac{-l}{2}}||\nabla\dot{\Delta}_lu||_{L^2}\leq2^\frac{-k}{2}||\nabla u||_{L^2},\\
||\nabla B^k||_{L^2}&\leq ||\nabla B||_{L^2},\\
||B^k||_{L^3}\leq\sum_{l\geq k}||\dot{\Delta}_lB||_{L^3}\leq\sum_{l\geq k} 2^{\frac{l}{2}}||\dot{\Delta}_lB||_{L^2}&\leq\sum_{l\geq k}2^{\frac{-l}{2}}||\nabla\dot{\Delta}_lB||_{L^2}\leq 2^{\frac{-k}{2}}||\nabla B||_{L^2},
\end{align*} here we have used the Bernstein inequalities \eqref{Ber1}-\eqref{Ber2} and H\"{o}lder inequality. From the embedding theorem
$||f||_{L^6(\mathbb{R}^3)}\leq C||\nabla f||_{L^2}$ and H\"{o}lder inequality, we also obtain the following estimates
\begin{align*}
||u\cdot\nabla u||_{L^\frac{3}{2}}&\leq C||u||_{L^6}||\nabla u||_{L^2}\leq C||\nabla u||^2_{L^2},\\
||B\cdot\nabla B||_{L^\frac{3}{2}}&\leq C||B||_{L^6}||\nabla B||_{L^2}\leq C||\nabla B||^2_{L^2},\\
||u\cdot\nabla B||_{L^\frac{3}{2}}&\leq C||u||_{L^6}||\nabla B||_{L^2}\leq C||\nabla u||_{L^2}||\nabla B||_{L^2},\\
||B\cdot\nabla u||_{L^\frac{3}{2}}&\leq C||B||_{L^6}||\nabla u||_{L^2}\leq C||\nabla u||_{L^2}||\nabla B||_{L^2}.
\end{align*}According to these obtained estimates, we can take $L^2$
inner product to the equations of velocity $u$ in \eqref{MHD} with $u^k$ then obtain that
\begin{align}\label{1}
\int_{\mathbb{R}^3}|\nabla u^k|^2dx=-\int_{\mathbb{R}^3}u\cdot\nabla u\cdot u^kdx+\int_{\mathbb{R}^3}B\cdot\nabla B\cdot u^kdx-\int_{\mathbb{R}^3}\nabla\dot{S}_ku:\nabla u^kdx.
\end{align}Similarly, Taking inner product to the equations of magnetic field $B$ in \eqref{MHD} with $B^k$
deduces that
\begin{align}\label{2}
\int_{\mathbb{R}^3}|\nabla B^k|^2dx=-\int_{\mathbb{R}^3}u\cdot\nabla B\cdot B^kdx+\int_{\mathbb{R}^3}B\cdot\nabla u\cdot B^kdx-\int_{\mathbb{R}^3}\nabla\dot{S}_kB:\nabla B^kdx.
\end{align}Adding \eqref{1} and \eqref{2}, we get that
\begin{align}\label{3}
&\int_{\mathbb{R}^3}(|\nabla u^k|^2+|\nabla B^k|^2)dx\\=&-\int_{\mathbb{R}^3}u\cdot\nabla u\cdot u^kdx-\int_{\mathbb{R}^3}u\cdot\nabla B\cdot B^kdx\notag
\\&+\int_{\mathbb{R}^3}B\cdot\nabla B\cdot u^kdx+\int_{\mathbb{R}^3}B\cdot\nabla u\cdot B^kdx\notag\\
&-\int_{\mathbb{R}^3}\nabla\dot{S}_ku:\nabla u^kdx-\int_{\mathbb{R}^3}\nabla\dot{S}_kB:\nabla B^kdx.\notag
\end{align}Noticing that $u=\dot{S}_ku+u^k$, $B=\dot{S}_kB+B^k$ and $\nabla\cdot u=\nabla\cdot B=0$, we deduce by using integration
by parts that
\begin{align}\label{4}
&-\int_{\mathbb{R}^3}u\cdot\nabla u\cdot u^kdx\\
=&-\int_{\mathbb{R}^3}\dot{S}_ku\cdot\nabla\dot{S}_ku\cdot u^kdx-\int_{\mathbb{R}^3}u^k\cdot\nabla\dot{S}_ku\cdot u^kdx-\int_{\mathbb{R}^3}u\cdot\nabla u^k\cdot u^kdx\notag\\
=&-\int_{\mathbb{R}^3}\dot{S}_ku\cdot\nabla\dot{S}_ku\cdot u^kdx-\int_{\mathbb{R}^3}u^k\cdot\nabla\dot{S}_ku\cdot u^kdx+\int_{\mathbb{R}^3}\nabla\cdot u\frac{1}{2}| u^k|^2dx\notag\\
=&-\int_{\mathbb{R}^3}\dot{S}_ku\cdot\nabla\dot{S}_ku\cdot u^kdx-\int_{\mathbb{R}^3}u^k\cdot\nabla\dot{S}_ku\cdot u^kdx,\notag
\end{align}
\begin{align}\label{5}
&-\int_{\mathbb{R}^3}u\cdot\nabla B\cdot B^kdx\\
=&-\int_{\mathbb{R}^3}\dot{S}_ku\cdot\nabla\dot{S}_kB\cdot B^kdx-\int_{\mathbb{R}^3}u^k\cdot\nabla\dot{S}_kB\cdot B^kdx-\int_{\mathbb{R}^3}u\cdot\nabla B^k\cdot B^kdx\notag\\
=&-\int_{\mathbb{R}^3}\dot{S}_ku\cdot\nabla\dot{S}_kB\cdot B^kdx-\int_{\mathbb{R}^3}u^k\cdot\nabla\dot{S}_kB\cdot B^kdx+\int_{\mathbb{R}^3}\nabla\cdot u\frac{1}{2}| B^k|^2dx\notag\\
=&-\int_{\mathbb{R}^3}\dot{S}_ku\cdot\nabla\dot{S}_kB\cdot B^kdx-\int_{\mathbb{R}^3}u^k\cdot\nabla\dot{S}_kB\cdot B^kdx,\notag
\end{align}
\begin{align}\label{6}
&\int_{\mathbb{R}^3}B\cdot\nabla B\cdot u^kdx+\int_{\mathbb{R}^3}B\cdot\nabla u\cdot B^kdx\\
&=\int_{\mathbb{R}^3}B\cdot\nabla B^k\cdot u^kdx+\int_{\mathbb{R}^3}B\cdot\nabla u^k\cdot B^kdx\notag\\&+\int_{\mathbb{R}^3}B\cdot\nabla\dot{S}_kB\cdot u^kdx
+\int_{\mathbb{R}^3}B\cdot\nabla\dot{S}_ku\cdot B^kdx\notag\\
&=-\int_{\mathbb{R}^3}B\cdot\nabla u^k\cdot B^kdx+\int_{\mathbb{R}^3}B\cdot\nabla u^k\cdot B^kdx\notag\\&+\int_{\mathbb{R}^3}B\cdot\nabla\dot{S}_kB\cdot u^kdx
+\int_{\mathbb{R}^3}B\cdot\nabla\dot{S}_ku\cdot B^kdx\notag\\
&=\int_{\mathbb{R}^3}\dot{S}_kB\cdot\nabla\dot{S}_kB\cdot u^kdx+\int_{\mathbb{R}^3}B^k\cdot\nabla\dot{S}_kB\cdot u^kdx\notag\\
&+\int_{\mathbb{R}^3}\dot{S}_kB\cdot\nabla\dot{S}_ku\cdot B^kdx+\int_{\mathbb{R}^3}B^k\cdot\nabla\dot{S}_ku\cdot B^kdx.\notag
\end{align}Substituting \eqref{4}-\eqref{6} into \eqref{3}, we get that
\begin{align}\label{7}
&\int_{\mathbb{R}^3}(|\nabla u^k|^2+|\nabla B^k|^2)dx\\
&=-\int_{\mathbb{R}^3}\dot{S}_ku\cdot\nabla\dot{S}_ku\cdot u^kdx-\int_{\mathbb{R}^3}u^k\cdot\nabla\dot{S}_ku\cdot u^kdx\notag\\
&-\int_{\mathbb{R}^3}\dot{S}_ku\cdot\nabla\dot{S}_kB\cdot B^kdx-\int_{\mathbb{R}^3}u^k\cdot\nabla\dot{S}_kB\cdot B^kdx\notag\\
&+\int_{\mathbb{R}^3}\dot{S}_kB\cdot\nabla\dot{S}_kB\cdot u^kdx+\int_{\mathbb{R}^3}B^k\cdot\nabla\dot{S}_kB\cdot u^kdx\notag\\
&+\int_{\mathbb{R}^3}\dot{S}_kB\cdot\nabla\dot{S}_ku\cdot B^kdx+\int_{\mathbb{R}^3}B^k\cdot\nabla\dot{S}_ku\cdot B^kdx\notag\\
&-\int_{\mathbb{R}^3}\nabla\dot{S}_ku:\nabla u^kdx-\int_{\mathbb{R}^3}\nabla\dot{S}_kB:\nabla B^kdx.\notag
\end{align}
It is not difficult to see that
\begin{align}
&\lim_{k\to-\infty}\int_{\mathbb{R}^3}(|\nabla u^k|^2+|\nabla B^k|^2)dx=\int_{\mathbb{R}^3}(|\nabla u|^2+|\nabla B|^2)dx\label{8},\\
&\lim_{k\to-\infty}\int_{\mathbb{R}^3}\nabla\dot{S}_ku:\nabla u^kdx=0\label{9},\\
&\lim_{k\to-\infty}\int_{\mathbb{R}^3}\nabla\dot{S}_kB:\nabla B^kdx=0\label{10}.
\end{align}Substituting \eqref{8}-\eqref{10} into \eqref{7}, it follows that
\begin{align}\label{11}
&\int_{\mathbb{R}^3}(|\nabla u|^2+|\nabla B|^2)dx\\=&\liminf_{k\to-\infty}\{-\int_{\mathbb{R}^3}\dot{S}_ku\cdot\nabla\dot{S}_ku\cdot u^kdx-\int_{\mathbb{R}^3}u^k\cdot\nabla\dot{S}_ku\cdot u^kdx\notag\\
&-\int_{\mathbb{R}^3}\dot{S}_ku\cdot\nabla\dot{S}_kB\cdot B^kdx-\int_{\mathbb{R}^3}u^k\cdot\nabla\dot{S}_kB\cdot B^kdx\notag\\
&+\int_{\mathbb{R}^3}\dot{S}_kB\cdot\nabla\dot{S}_kB\cdot u^kdx+\int_{\mathbb{R}^3}B^k\cdot\nabla\dot{S}_kB\cdot u^kdx\notag\\
&+\int_{\mathbb{R}^3}\dot{S}_kB\cdot\nabla\dot{S}_ku\cdot B^kdx+\int_{\mathbb{R}^3}B^k\cdot\nabla\dot{S}_ku\cdot B^kdx,\notag\}\\
=&\liminf_{k\to-\infty}(I_1+I_2+I_3+I_4+I_5+I_6+I_7+I_8).\notag
\end{align}

Next, we will estimate the terms $I_1-I_8$. Throughout, we will use Einstein summation convention (summing over repeated indices).

For $I_1$, by Bony decomposition, we get
\begin{align}\label{12}
&-I_1\\=&\int_{\mathbb{R}^3}\partial_j\dot{S}_ku_i\dot{S}_ku_j u^k_idx\notag\\
=&\int_{\mathbb{R}^3}\partial_j\dot{S}_ku_i(\sum_{l\in\mathbb{Z}}\dot{\Delta}_l\dot{S}_ku_j\dot{S}_{l-2}u^k_i+\sum_{l\in\mathbb{Z}}\dot{\Delta}_lu^k_i\dot{S}_{l-2}\dot{S}_ku_j
+\sum_{l\in\mathbb{Z}}\dot{\Delta}_l\dot{S}_ku_j\tilde{\dot{\Delta}}_lu^k_i)dx\notag\\
=&I_{11}+I_{12}+I_{13},\notag
\end{align} where we have used the notation $\tilde{\dot{\Delta}}_lf=\sum_{|l'-l|\leq 2}\dot{\Delta}_{l'}f$.

For $I_{11}$, we observe that $\dot{\Delta}_l\dot{S}_ku_j\dot{S}_{l-2}u^k_i\neq0$ means $l\leq k$ and $l\geq k+2$. This means that $\dot{\Delta}_l\dot{S}_ku_j\dot{S}_{l-2}u^k_i=0$ for all $l\in\mathbb{Z}$ and then
\begin{align*}
I_{11}=\int_{\mathbb{R}^3}\partial_j\dot{S}_ku_i\sum_{l\in\mathbb{Z}}\dot{\Delta}_l\dot{S}_ku_j\dot{S}_{l-2}u^k_idx=0.
\end{align*}For $I_{12}$, we first observe that $\dot{\Delta}_lu^k_i\dot{S}_{l-2}\dot{S}_ku_j\neq0$ means $l\geq k-1$ and
\begin{align}\label{13}
supp \mathcal{F}(\dot{\Delta}_lu^k_i\dot{S}_{l-2}\dot{S}_ku_j)\subset\{\xi:2^{l-2}\leq|\xi|<\frac{9}{8}2^{l+1}\}.
\end{align} On the other hand, we also have
\begin{align}\label{14}
supp \mathcal{F}(\partial_j\dot{S}_ku_i)\subset\{\xi:|\xi|< 2^k\}.
\end{align}Combine \eqref{13} and \eqref{14} implies that
\begin{align*}
&I_{12}\\=&\int_{\mathbb{R}^3}\partial_j\dot{S}_ku_i\sum_{l\in\mathbb{Z}}\dot{\Delta}_lu^k_i\dot{S}_{l-2}\dot{S}_ku_jdx
=\int_{\mathbb{R}^3}\partial_j\dot{S}_ku_i\sum_{l=k-1}^{k+1}\dot{\Delta}_lu^k_i\dot{S}_{l-2}\dot{S}_{k}u_jdx\\
=&\int_{\mathbb{R}^3}\sum_{l=k-1}^{k+1}\sum_{l'=l-2}^{k-1}\partial_j\dot{\Delta}_{l'}u_i\dot{\Delta}_lu^k_i\dot{S}_{l-2}\dot{S}_{k}u_jdx\notag\\=&\int_{\mathbb{R}^3}
\sum_{l=k-1}^{k+1}\sum_{l'=l-2}^{k-1}\dot{S}_{l-2}\dot{S}_{k}u\cdot\nabla\dot{\Delta}_{l'}u\cdot\dot{\Delta}_lu^kdx.\notag
\end{align*} We now investigate $I_{13}$. It is not difficult to see that $\dot{\Delta}_l\dot{S}_ku_j\neq0\Rightarrow l\leq k$ and $\tilde{\dot{\Delta}}_lu^k_i\neq0\Rightarrow l\geq k-3$.
Hence it follows from the above facts that
\begin{align*}
\sum_{l\in\mathbb{Z}}\dot{\Delta}_l\dot{S}_ku_j\tilde{\dot{\Delta}}_lu^k_i=\sum_{l=k-3}^{k}\dot{\Delta}_l\dot{S}_ku_j\tilde{\dot{\Delta}}_lu^k_i.
\end{align*}From the above identity, we deduce that
\begin{align*}
I_{13}=&\int_{\mathbb{R}^3}\partial_j\dot{S}_ku_i\sum_{l\in\mathbb{Z}}\dot{\Delta}_l\dot{S}_ku_j\tilde{\dot{\Delta}}_lu^k_i
=\int_{\mathbb{R}^3}\partial_j\dot{S}_ku_i\sum_{l=k-3}^{k}\dot{\Delta}_l\dot{S}_ku_j\tilde{\dot{\Delta}}_lu^k_idx\\
=&\int_{\mathbb{R}^3}\sum_{l=k-3}^{k}\dot{\Delta}_l\dot{S}_ku\cdot\nabla\dot{S}_ku\cdot\tilde{\dot{\Delta}}_lu^kdx.\notag
\end{align*}Substituting the estimates of $I_{11}$, $I_{12}$ and $I_{13}$ into \eqref{12}, it follows that
\begin{align}\label{15}
&-I_1\\=&\int_{\mathbb{R}^3}
\sum_{l=k-1}^{k+1}\sum_{l'=l-2}^{k-1}\dot{S}_{l-2}\dot{S}_{k}u\cdot\nabla\dot{\Delta}_{l'}u\cdot\dot{\Delta}_lu^kdx
+\int_{\mathbb{R}^3}\sum_{l=k-3}^{k}\dot{\Delta}_l\dot{S}_ku\cdot\nabla\dot{S}_ku\cdot\tilde{\dot{\Delta}}_lu^kdx.\notag
\end{align}

We now restrict attention to the term $I_2$. Similarly,  by using Bony decomposition, it follows that
\begin{align}\label{16}
&-I_2\\=&\int_{\mathbb{R}^3}\partial_j\dot{S}_ku_iu^k_ju^k_idx\notag\\
=&\int_{\mathbb{R}^3}\partial_j\dot{S}_ku_i(\sum_{l\in\mathbb{Z}}\dot{\Delta}_lu^k_j\dot{S}_{l-2}u^k_i+\sum_{l\in\mathbb{Z}}\dot{\Delta}_lu^k_i\dot{S}_{l-2}u^k_j
+\sum_{l\in\mathbb{Z}}\dot{\Delta}_lu^k_i\tilde{\dot{\Delta}}_{l}u^k_j)dx\notag\\
=&I_{21}+I_{22}+I_{23}.\notag
\end{align}
For $I_{21}$, we first observe that $\dot{\Delta}_lu^k_j\dot{S}_{l-2}u^k_i\neq0\Rightarrow l\geq k+2$ and then
\begin{align}\label{17}
supp\mathcal{F}(\dot{\Delta}_lu^k_j\dot{S}_{l-2}u^k_i)\subset\{\xi:2^{l-2}\leq|\xi|<\frac{9}{8}2^{l+1}\} \quad for\quad l\geq k+2.
\end{align}Combining \eqref{14} and \eqref{17} shows that
\begin{align}\label{I21}
I_{21}=&\int_{\mathbb{R}^3}\partial_j\dot{S}_ku_i\sum_{l\in\mathbb{Z}}\dot{\Delta}_lu^k_j\dot{S}_{l-2}u^k_idx\\
=&\int_{\mathbb{R}^3}\partial_j\dot{S}_ku_i\sum_{l\geq k+2}\dot{\Delta}_lu^k_j\dot{S}_{l-2}u^k_idx\notag\\
=&\int_{\mathbb{R}^3}\mathcal{F}(\partial_j\dot{S}_ku_i)\mathcal{F}(\sum_{l\geq k+2}\dot{\Delta}_lu^k_j\dot{S}_{l-2}u^k_i)d\xi\notag\\
=&0.\notag
\end{align}Applying the same arguments, we also obtain
\begin{align}\label{I22}
I_{22}=\int_{\mathbb{R}^3}\partial_j\dot{S}_ku_i\sum_{l\in\mathbb{Z}}\dot{\Delta}_lu^k_i\dot{S}_{l-2}u^k_jdx=0.
\end{align}

Finally, we consider the most complicated quantity $I_{23}$. It is not difficult to see that $\dot{\Delta}_lu^k_i\neq0\Rightarrow l\geq k-1$, we thus get that
\begin{align*}
I_{23}=\int_{\mathbb{R}^3}\partial_j\dot{S}_ku_i\sum_{l\in\mathbb{Z}}\dot{\Delta}_lu^k_i\tilde{\dot{\Delta}}_{l}u^k_jdx
=\int_{\mathbb{R}^3}\partial_j\dot{S}_ku_i\sum_{l\geq k-1}\dot{\Delta}_lu^k_i\tilde{\dot{\Delta}}_{l}u^k_jdx.
\end{align*} Since $supp\mathcal{F}(\dot{\Delta}_lu^k_i\tilde{\dot{\Delta}}_{l}u^k_j)\subset\{\xi:|\xi|<5\times2^{l+1}\}$, we can not localize $I_{23}$ into the region near the origin in frequency space directly.
The key observation is that we can decompose $I_{23}$ into two parts, the first one is a low frequency part compared to $k$, the other is a high frequency part compared to $k$ and the high frequency part will vanish as $k\to-\infty$. Let $\theta\in (0,1)$ whose value will be determined later. We decompose $I_{23}$ as follows
\begin{align}\label{deI23}
I_{23}=&\int_{\mathbb{R}^3}\partial_j\dot{S}_ku_i\sum_{k-1\leq l\leq [\theta k]}\dot{\Delta}_lu^k_i\tilde{\dot{\Delta}}_{l}u^k_jdx+\int_{\mathbb{R}^3}\partial_j\dot{S}_ku_i\sum_{l\geq [\theta k]+1}\dot{\Delta}_lu^k_i\tilde{\dot{\Delta}}_{l}u^k_jdx\\
=&I_{231}+I_{232}.\notag
\end{align} The quantity $I_{232}$ can be bounded as following:
\begin{align}\label{bI232}
I_{232}&\leq||\partial_j\dot{S}_ku_i||_{L^\infty}\sum_{l\geq [\theta k]+1}||\dot{\Delta}_lu^k_i||_{L^2}||\tilde{\dot{\Delta}}_{l}u^k_j||_{L^2}\\
&\leq 2^{\frac{3}{2}k}||\nabla\dot{S}_ku_i||_{L^2}\sum_{l\geq [\theta k]+1}2^{-2l}||\nabla\dot{\Delta}_lu^k_i||_{L^2}||\nabla\tilde{\dot{\Delta}}_{l}u^k_j||_{L^2}\notag\\
&\leq 2^{(\frac{3}{2}-2\theta)k}||\nabla\dot{S}_ku_i||_{L^2}\sum_{l\geq [\theta k]+1}2^{2(\theta k-l)}||\nabla\dot{\Delta}_lu^k_i||_{L^2}||\nabla\tilde{\dot{\Delta}}_{l}u^k_j||_{L^2}\notag\\
&\leq C2^{(\frac{3}{2}-2\theta )k}||\nabla\dot{S}_ku_i||_{L^2}||\nabla u||^2_{L^2},\notag
\end{align} where we have used \eqref{Ber1}-\eqref{Ber2}. Since $||\nabla\dot{S}_ku_i||_{L^2}\to0$, as
$k\to-\infty$, we deduce that
\begin{align*}
\lim_{k\to-\infty}I_{232}=0 \quad if \quad\theta\leq\frac{3}{4}.
\end{align*}Taking $\theta=\frac{3}{4}$, substituting the estimates of $I_{21}$, $I_{22}$, $I_{23}$ and $I_{232}$ into \eqref{16} implies that
\begin{align}\label{I_2}
&-I_2\\=&\int_{\mathbb{R}^3}\partial_j\dot{S}_ku_i\sum_{l=k-1}^{[\frac{3}{4} k]}\dot{\Delta}_lu^k_i\tilde{\dot{\Delta}}_{l}u^k_jdx+\int_{\mathbb{R}^3}\partial_j\dot{S}_ku_i
\sum_{l\geq [\frac{3}{4} k]+1}\dot{\Delta}_lu^k_i\tilde{\dot{\Delta}}_{l}u^k_jdx\notag\\
=&\int_{\mathbb{R}^3}\sum_{l=k-1}^{[\frac{3}{4} k]}\tilde{\dot{\Delta}}_{l}u^k\cdot\nabla\dot{S}_ku\cdot\dot{\Delta}_lu^kdx+\int_{\mathbb{R}^3}
\sum_{l\geq [\frac{3}{4} k]+1}\tilde{\dot{\Delta}}_{l}u^k\cdot\nabla\dot{S}_ku\cdot\dot{\Delta}_lu^kdx\notag
\end{align} and
\begin{align}\label{vI_2}
\lim_{k\to-\infty}|\int_{\mathbb{R}^3}
\sum_{l\geq [\frac{3}{4} k]+1}\tilde{\dot{\Delta}}_{l}u^k\cdot\nabla\dot{S}_ku\cdot\dot{\Delta}_lu^kdx|=0.
\end{align}

 Observe that the terms $I_1,I_3,I_5,I_7$ are with the same structure $\int_{\mathbb{R}^3}\dot{S}_k f\cdot\nabla\dot{S}_kg\cdot h^kdx$. Applying the same arguments of $I_1$ to the terms $I_3$, $I_5$ and $I_7$, we obtain that
\begin{align}\label{I_3}
&-I_3\\=&\int_{\mathbb{R}^3}
\sum_{l=k-1}^{k+1}\sum_{l'=l-2}^{k-1}\dot{S}_{l-2}\dot{S}_{k}u\cdot\nabla\dot{\Delta}_{l'}B\cdot\dot{\Delta}_lB^kdx
+\int_{\mathbb{R}^3}\sum_{l=k-3}^{k}\dot{\Delta}_l\dot{S}_ku\cdot\nabla\dot{S}_kB\cdot\tilde{\dot{\Delta}}_lB^kdx.\notag
\end{align}
\begin{align}\label{I_5}
&I_5\\=&\int_{\mathbb{R}^3}
\sum_{l=k-1}^{k+1}\sum_{l'=l-2}^{k-1}\dot{S}_{l-2}\dot{S}_{k}B\cdot\nabla\dot{\Delta}_{l'}B\cdot\dot{\Delta}_lu^kdx
+\int_{\mathbb{R}^3}\sum_{l=k-3}^{k}\dot{\Delta}_l\dot{S}_kB\cdot\nabla\dot{S}_kB\cdot\tilde{\dot{\Delta}}_lu^kdx.\notag
\end{align}
\begin{align}\label{I_7}
&I_7\\=&\int_{\mathbb{R}^3}
\sum_{l=k-1}^{k+1}\sum_{l'=l-2}^{k-1}\dot{S}_{l-2}\dot{S}_{k}B\cdot\nabla\dot{\Delta}_{l'}u\cdot\dot{\Delta}_lB^kdx
+\int_{\mathbb{R}^3}\sum_{l=k-3}^{k}\dot{\Delta}_l\dot{S}_kB\cdot\nabla\dot{S}_ku\cdot\tilde{\dot{\Delta}}_lB^kdx.\notag
\end{align}

We also observe that the terms $I_2,I_4,I_6,I_8$ are with the same structure $\int_{\mathbb{R}^3}f^k\cdot\nabla\dot{S}_kg\cdot h^kdx$. Applying the same arguments of $I_2$ to the terms $I_4$, $I_6$ and $I_8$, we obtain that
\begin{align}\label{I_4}
&-I_4\\=&\int_{\mathbb{R}^3}\sum_{l=k-1}^{[\frac{3}{4} k]}\tilde{\dot{\Delta}}_{l}u^k\cdot\nabla\dot{S}_kB\cdot\dot{\Delta}_lB^kdx+\int_{\mathbb{R}^3}
\sum_{l\geq [\frac{3}{4} k]+1}\tilde{\dot{\Delta}}_{l}u^k\cdot\nabla\dot{S}_kB\cdot\dot{\Delta}_lB^kdx\notag
\end{align} and
\begin{align}\label{vI_4}
\lim_{k\to-\infty}|\int_{\mathbb{R}^3}
\sum_{l\geq [\frac{3}{4} k]+1}\tilde{\dot{\Delta}}_{l}u^k\cdot\nabla\dot{S}_kB\cdot\dot{\Delta}_lB^kdx|=0;
\end{align}
\begin{align}\label{I_6}
&I_6\\=&\int_{\mathbb{R}^3}\sum_{l=k-1}^{[\frac{3}{4} k]}\tilde{\dot{\Delta}}_{l}B^k\cdot\nabla\dot{S}_kB\cdot\dot{\Delta}_lu^kdx+\int_{\mathbb{R}^3}
\sum_{l\geq [\frac{3}{4} k]+1}\tilde{\dot{\Delta}}_{l}B^k\cdot\nabla\dot{S}_kB\cdot\dot{\Delta}_lu^kdx\notag
\end{align} and
\begin{align}\label{vI_6}
\lim_{k\to-\infty}|\int_{\mathbb{R}^3}
\sum_{l\geq [\frac{3}{4} k]+1}\tilde{\dot{\Delta}}_{l}B^k\cdot\nabla\dot{S}_kB\cdot\dot{\Delta}_lu^kdx|=0;
\end{align}
\begin{align}\label{I_8}
&I_8\\=&\int_{\mathbb{R}^3}\sum_{l=k-1}^{[\frac{3}{4} k]}\tilde{\dot{\Delta}}_{l}B^k\cdot\nabla\dot{S}_ku\cdot\dot{\Delta}_lB^kdx+\int_{\mathbb{R}^3}
\sum_{l\geq [\frac{3}{4} k]+1}\tilde{\dot{\Delta}}_{l}B^k\cdot\nabla\dot{S}_ku\cdot\dot{\Delta}_lB^kdx\notag
\end{align} and
\begin{align}\label{vI_8}
\lim_{k\to-\infty}|\int_{\mathbb{R}^3}
\sum_{l\geq [\frac{3}{4} k]+1}\tilde{\dot{\Delta}}_{l}B^k\cdot\nabla\dot{S}_ku\cdot\dot{\Delta}_lB^kdx|=0.
\end{align}

Substituting \eqref{15} and \eqref{I_2}-\eqref{vI_8} into \eqref{11}, we conclude that
\begin{align*}
&\int_{\mathbb{R}^3}(|\nabla u|^2+|\nabla B|^2)dx\\=&\liminf_{k\to-\infty}\{-(\int_{\mathbb{R}^3}
\sum_{l=k-1}^{k+1}\sum_{l'=l-2}^{k-1}\dot{S}_{l-2}\dot{S}_{k}u\cdot\nabla\dot{\Delta}_{l'}u\cdot\dot{\Delta}_lu^kdx
+\int_{\mathbb{R}^3}\sum_{l=k-3}^{k}\dot{\Delta}_l\dot{S}_ku\cdot\nabla\dot{S}_ku\cdot\tilde{\dot{\Delta}}_lu^kdx)\notag\\
&-(\int_{\mathbb{R}^3}\sum_{l=k-1}^{k+1}\sum_{l'=l-2}^{k-1}\dot{S}_{l-2}\dot{S}_{k}u\cdot\nabla\dot{\Delta}_{l'}B\cdot\dot{\Delta}_lB^kdx
+\int_{\mathbb{R}^3}\sum_{l=k-3}^{k}\dot{\Delta}_l\dot{S}_ku\cdot\nabla\dot{S}_kB\cdot\tilde{\dot{\Delta}}_lB^kdx)\notag\\
+&(\int_{\mathbb{R}^3}\sum_{l=k-1}^{k+1}\sum_{l'=l-2}^{k-1}\dot{S}_{l-2}\dot{S}_{k}B\cdot\nabla\dot{\Delta}_{l'}B\cdot\dot{\Delta}_lu^kdx
+\int_{\mathbb{R}^3}\sum_{l=k-3}^{k}\dot{\Delta}_l\dot{S}_kB\cdot\nabla\dot{S}_kB\cdot\tilde{\dot{\Delta}}_lu^kdx)\notag\\
&+(\int_{\mathbb{R}^3}\sum_{l=k-1}^{k+1}\sum_{l'=l-2}^{k-1}\dot{S}_{l-2}\dot{S}_{k}B\cdot\nabla\dot{\Delta}_{l'}u\cdot\dot{\Delta}_lB^kdx
+\int_{\mathbb{R}^3}\sum_{l=k-3}^{k}\dot{\Delta}_l\dot{S}_kB\cdot\nabla\dot{S}_ku\cdot\tilde{\dot{\Delta}}_lB^kdx)\notag\\
&-\int_{\mathbb{R}^3}\sum_{l=k-1}^{[\frac{3}{4} k]}\tilde{\dot{\Delta}}_{l}u^k\cdot\nabla\dot{S}_ku\cdot\dot{\Delta}_lu^kdx-\int_{\mathbb{R}^3}\sum_{l=k-1}^{[\frac{3}{4} k]}\tilde{\dot{\Delta}}_{l}u^k\cdot\nabla\dot{S}_kB\cdot\dot{\Delta}_lB^kdx\notag\\
&+\int_{\mathbb{R}^3}\sum_{l=k-1}^{[\frac{3}{4} k]}\tilde{\dot{\Delta}}_{l}B^k\cdot\nabla\dot{S}_kB\cdot\dot{\Delta}_lu^kdx+\int_{\mathbb{R}^3}\sum_{l=k-1}^{[\frac{3}{4} k]}\tilde{\dot{\Delta}}_{l}B^k\cdot\nabla\dot{S}_ku\cdot\dot{\Delta}_lB^kdx\notag\}.
\end{align*} This completes the proof of Lemma \ref{lemma}.
\qquad $\hfill\Box$

Next, we will prove the first Liouville type theorem based on Lemma \ref{lemma}.

\emph{\bf{Proof of Theorem \ref{main}.}} From Lemma \ref{lemma}, we know that
\begin{align}\label{+11}
\int_{\mathbb{R}^3}(|\nabla u|^2+|\nabla B|^2)dx=\liminf_{k\to-\infty}\sum_{i=1}^{8}J_i
\end{align}with
\begin{align*}
J_1=&-(\int_{\mathbb{R}^3}
\sum_{l=k-1}^{k+1}\sum_{l'=l-2}^{k-1}\dot{S}_{l-2}\dot{S}_{k}u\cdot\nabla\dot{\Delta}_{l'}u\cdot\dot{\Delta}_lu^kdx
+\int_{\mathbb{R}^3}\sum_{l=k-3}^{k}\dot{\Delta}_l\dot{S}_ku\cdot\nabla\dot{S}_ku\cdot\tilde{\dot{\Delta}}_lu^kdx)\\
J_2=&-(\int_{\mathbb{R}^3}\sum_{l=k-1}^{k+1}\sum_{l'=l-2}^{k-1}\dot{S}_{l-2}\dot{S}_{k}u\cdot\nabla\dot{\Delta}_{l'}B\cdot\dot{\Delta}_lB^kdx
+\int_{\mathbb{R}^3}\sum_{l=k-3}^{k}\dot{\Delta}_l\dot{S}_ku\cdot\nabla\dot{S}_kB\cdot\tilde{\dot{\Delta}}_lB^kdx)\\
J_3=&\int_{\mathbb{R}^3}\sum_{l=k-1}^{k+1}\sum_{l'=l-2}^{k-1}\dot{S}_{l-2}\dot{S}_{k}B\cdot\nabla\dot{\Delta}_{l'}B\cdot\dot{\Delta}_lu^kdx
+\int_{\mathbb{R}^3}\sum_{l=k-3}^{k}\dot{\Delta}_l\dot{S}_kB\cdot\nabla\dot{S}_kB\cdot\tilde{\dot{\Delta}}_lu^kdx\\
J_4=&\int_{\mathbb{R}^3}\sum_{l=k-1}^{k+1}\sum_{l'=l-2}^{k-1}\dot{S}_{l-2}\dot{S}_{k}B\cdot\nabla\dot{\Delta}_{l'}u\cdot\dot{\Delta}_lB^kdx
+\int_{\mathbb{R}^3}\sum_{l=k-3}^{k}\dot{\Delta}_l\dot{S}_kB\cdot\nabla\dot{S}_ku\cdot\tilde{\dot{\Delta}}_lB^kdx
\end{align*}
\begin{align*}
J_5=&-\int_{\mathbb{R}^3}\sum_{l=k-1}^{[\frac{3}{4} k]}\tilde{\dot{\Delta}}_{l}u^k\cdot\nabla\dot{S}_ku\cdot\dot{\Delta}_lu^kdx\quad\quad
J_6=-\int_{\mathbb{R}^3}\sum_{l=k-1}^{[\frac{3}{4} k]}\tilde{\dot{\Delta}}_{l}u^k\cdot\nabla\dot{S}_kB\cdot\dot{\Delta}_lB^kdx\\
J_7=&\int_{\mathbb{R}^3}\sum_{l=k-1}^{[\frac{3}{4} k]}\tilde{\dot{\Delta}}_{l}B^k\cdot\nabla\dot{S}_kB\cdot\dot{\Delta}_lu^kdx\quad\quad
J_8=\int_{\mathbb{R}^3}\sum_{l=k-1}^{[\frac{3}{4} k]}\tilde{\dot{\Delta}}_{l}B^k\cdot\nabla\dot{S}_ku\cdot\dot{\Delta}_lB^kdx.
\end{align*}
Now, our aims are to bound the right hand side terms in \eqref{+11} based on Bernstein inequalities \eqref{Ber1}-\eqref{Ber2} and \eqref{+10}.
For the terms $J_1-J_8$, we estimate them as follows:
\begin{align}\label{b1}
&|J_1|\\=&|\int_{\mathbb{R}^3}
\sum_{l=k-1}^{k+1}\sum_{l'=l-2}^{k-1}\dot{S}_{l-2}\dot{S}_{k}u\cdot\nabla\dot{\Delta}_{l'}u\cdot\dot{\Delta}_lu^kdx+\int_{\mathbb{R}^3}\sum_{l=k-3}^{k}
\dot{\Delta}_l\dot{S}_ku\cdot\nabla\dot{S}_ku\cdot\tilde{\dot{\Delta}}_lu^kdx|\notag\\
\leq&\sum_{l=k-1}^{k+1}\sum_{l'=l-2}^{k-1}||\dot{S}_{l-2}\dot{S}_ku||_{L^\infty}||\nabla\dot{\Delta}_{l'}u||_{L^2}||\dot{\Delta}_lu^k||_{L^2}
+\sum_{l=k-3}^{k}||\nabla\dot{S}_{k}u||_{L^\infty}||\dot{\Delta}_{l}\dot{S}_ku||_{L^2}||\tilde{\dot{\Delta}}_lu^k||_{L^2}\notag\\
\leq&\sum_{l=k-1}^{k+1}\sum_{l'=l-2}^{k-1}||\dot{S}_ku||_{L^\infty}||\nabla\dot{\Delta}_{l'}u||_{L^2}||\dot{\Delta}_lu||_{L^2}
+\sum_{l=k-3}^{k}2^k||\dot{S}_{k}u||_{L^\infty}||\dot{\Delta}_{l}u||_{L^2}||\tilde{\dot{\Delta}}_lu||_{L^2}\notag\\
\leq&2^{-k}||\dot{S}_ku||_{L^\infty}(\sum_{l=k-1}^{k+1}\sum_{l'=l-2}^{k-1}2^{k-l}||\nabla\dot{\Delta}_{l'}u||_{L^2}||\nabla\dot{\Delta}_lu||_{L^2}
+\sum_{l=k-3}^{k}2^{2k-2l}||\nabla\dot{\Delta}_{l}u||_{L^2}||\nabla\tilde{\dot{\Delta}}_lu||_{L^2})\notag\\
\leq& C2^{-k}||\dot{S}_{k}u||_{L^\infty}||\nabla\dot{S}_{k+3}u||^2_{L^2},\notag
\end{align}
\begin{align}\label{b2}
&|J_2|\\=&\int_{\mathbb{R}^3}\sum_{l=k-1}^{k+1}\sum_{l'=l-2}^{k-1}\dot{S}_{l-2}\dot{S}_{k}u\cdot\nabla\dot{\Delta}_{l'}B\cdot\dot{\Delta}_lB^kdx
+\int_{\mathbb{R}^3}\sum_{l=k-3}^{k}\dot{\Delta}_l\dot{S}_ku\cdot\nabla\dot{S}_kB\cdot\tilde{\dot{\Delta}}_lB^kdx|\notag\\
\leq&\sum_{l=k-1}^{k+1}\sum_{l'=l-2}^{k-1}||\dot{S}_{l-2}\dot{S}_ku||_{L^\infty}||\nabla\dot{\Delta}_{l'}B||_{L^2}||\dot{\Delta}_lB^k||_{L^2}
+\sum_{l=k-3}^{k}||\nabla\dot{S}_{k}B||_{L^\infty}||\dot{\Delta}_{l}\dot{S}_ku||_{L^2}||\tilde{\dot{\Delta}}_lB^k||_{L^2}\notag\\
\leq&\sum_{l=k-1}^{k+1}\sum_{l'=l-2}^{k-1}||\dot{S}_ku||_{L^\infty}||\nabla\dot{\Delta}_{l'}B||_{L^2}||\dot{\Delta}_lB||_{L^2}
+\sum_{l=k-3}^{k}2^k||\dot{S}_{k}B||_{L^\infty}||\dot{\Delta}_{l}u||_{L^2}||\tilde{\dot{\Delta}}_lB||_{L^2}\notag\\
\leq& C2^{-k}(||\dot{S}_{k}u||_{L^\infty}+||\dot{S}_{k}B||_{L^\infty})(||\nabla\dot{S}_{k+3}B||^2_{L^2}+||\nabla\dot{S}_{k+3}u||^2_{L^2}),\notag
\end{align}
\begin{align}\label{b3}
&|J_3|\\=&
|\int_{\mathbb{R}^3}\sum_{l=k-1}^{k+1}\sum_{l'=l-2}^{k-1}\dot{S}_{l-2}\dot{S}_{k}B\cdot\nabla\dot{\Delta}_{l'}B\cdot\dot{\Delta}_lu^kdx
+\int_{\mathbb{R}^3}\sum_{l=k-3}^{k}\dot{\Delta}_l\dot{S}_kB\cdot\nabla\dot{S}_kB\cdot\tilde{\dot{\Delta}}_lu^kdx|\notag\\
\leq&\sum_{l=k-1}^{k+1}\sum_{l'=l-2}^{k-1}||\dot{S}_{l-2}\dot{S}_kB||_{L^\infty}||\nabla\dot{\Delta}_{l'}B||_{L^2}||\dot{\Delta}_lu^k||_{L^2}
+\sum_{l=k-3}^{k}||\nabla\dot{S}_{k}B||_{L^\infty}||\dot{\Delta}_{l}\dot{S}_kB||_{L^2}||\tilde{\dot{\Delta}}_lu^k||_{L^2}\notag\\
\leq&\sum_{l=k-1}^{k+1}\sum_{l'=l-2}^{k-1}||\dot{S}_kB||_{L^\infty}||\nabla\dot{\Delta}_{l'}B||_{L^2}||\dot{\Delta}_lu||_{L^2}
+C\sum_{l=k-3}^{k}2^k||\dot{S}_{k}B||_{L^\infty}||\dot{\Delta}_{l}B||_{L^2}||\tilde{\dot{\Delta}}_lu||_{L^2}\notag\\
\leq&2^{-k}||\dot{S}_kB||_{L^\infty}(\sum_{l=k-1}^{k+1}\sum_{l'=l-2}^{k-1}||\nabla\dot{\Delta}_{l'}B||_{L^2}||\nabla\dot{\Delta}_lu||_{L^2}
+\sum_{l=k-3}^{k}||\nabla\dot{\Delta}_{l}B||_{L^2}||\nabla\tilde{\dot{\Delta}}_lu||_{L^2})\notag\\
\leq& C2^{-k}||\dot{S}_{k}B||_{L^\infty}(||\nabla\dot{S}_{k+3}B||^2_{L^2}+||\nabla\dot{S}_{k+3}u||^2_{L^2}),\notag
\end{align}
\begin{align}\label{b4}
&|J_4|\\=&
|\int_{\mathbb{R}^3}\sum_{l=k-1}^{k+1}\sum_{l'=l-2}^{k-1}\dot{S}_{l-2}\dot{S}_{k}B\cdot\nabla\dot{\Delta}_{l'}u\cdot\dot{\Delta}_lB^kdx
+\int_{\mathbb{R}^3}\sum_{l=k-3}^{k}\dot{\Delta}_l\dot{S}_kB\cdot\nabla\dot{S}_ku\cdot\tilde{\dot{\Delta}}_lB^kdx|\notag\\
\leq&\sum_{l=k-1}^{k+1}\sum_{l'=l-2}^{k-1}||\dot{S}_{l-2}\dot{S}_kB||_{L^\infty}||\nabla\dot{\Delta}_{l'}u||_{L^2}||\dot{\Delta}_lB^k||_{L^2}
+\sum_{l=k-3}^{k}||\nabla\dot{S}_{k}u||_{L^\infty}||\dot{\Delta}_{l}\dot{S}_kB||_{L^2}||\tilde{\dot{\Delta}}_lB^k||_{L^2}\notag\\
\leq&C\sum_{l=k-1}^{k+1}\sum_{l'=l-2}^{k-1}||\dot{S}_kB||_{L^\infty}||\nabla\dot{\Delta}_{l'}u||_{L^2}||\dot{\Delta}_lB||_{L^2}
+C\sum_{l=k-3}^{k}2^k||\dot{S}_{k}u||_{L^\infty}||\dot{\Delta}_{l}B||_{L^2}||\tilde{\dot{\Delta}}_lB||_{L^2}\notag\\
\leq& C2^{-k}(||\dot{S}_{k}u||_{L^\infty}+||\dot{S}_{k}B||_{L^\infty})(||\nabla\dot{S}_{k+3}B||^2_{L^2}+||\nabla\dot{S}_{k+3}u||^2_{L^2}).\notag
\end{align}

\begin{align}\label{b5}
|J_5|=&|\int_{\mathbb{R}^3}\sum_{l=k-1}^{[\frac{3k}{4}]}\dot{\Delta}_l u^k\cdot\nabla \dot{S}_ku\cdot\tilde{\dot{\Delta}}_lu^kdx|\\
\leq&||\nabla\dot{S}_ku||_{L^\infty}\sum_{l=k-1}^{[\frac{3k}{4}]}||\dot{\Delta}_l u||_{L^2}||\tilde{\dot{\Delta}}_lu||_{L^2}\notag\\
\leq& C2^k||\dot{S}_ku||_{L^\infty}\sum_{l=k-1}^{[\frac{3k}{4}]}2^{-2l}||\nabla\dot{\Delta}_l u||_{L^2}||\nabla\tilde{\dot{\Delta}}_lu||_{L^2}\notag\\
\leq& C2^{-k}||\dot{S}_ku||_{L^\infty}\sum_{l=k-1}^{[\frac{3k}{4}]}2^{2(k-l)}||\nabla\dot{\Delta}_l u||_{L^2}||\nabla\tilde{\dot{\Delta}}_lu||_{L^2}\notag\\
\leq&C2^{-k}||\dot{S}_ku||_{L^\infty}||\nabla\dot{S}_{[\frac{3k}{4}]+3}u||^2_{L^2},\notag
\end{align}
\begin{align}\label{b6}
|J_6|=&|\int_{\mathbb{R}^3}\sum_{l=k-1}^{[\frac{3}{4} k]}\tilde{\dot{\Delta}}_{l}u^k\cdot\nabla\dot{S}_kB\cdot\dot{\Delta}_lB^kdx|\\
\leq&||\nabla\dot{S}_kB||_{L^\infty}\sum_{l=k-1}^{[\frac{3k}{4}]}||\dot{\Delta}_l B^k||_{L^2}||\tilde{\dot{\Delta}}_lu^k||_{L^2}\notag\\
\leq& C2^k||\dot{S}_kB||_{L^\infty}\sum_{l=k-1}^{[\frac{3k}{4}]}2^{-2l}||\nabla\dot{\Delta}_l B||_{L^2}||\nabla\tilde{\dot{\Delta}}_lu||_{L^2}\notag\\
\leq& C2^{-k}||\dot{S}_kB||_{L^\infty}\sum_{l=k-1}^{[\frac{3k}{4}]}2^{2(k-l)}||\nabla\dot{\Delta}_l B||_{L^2}||\nabla\tilde{\dot{\Delta}}_lu||_{L^2}\notag\\
\leq&C2^{-k}||\dot{S}_kB||_{L^\infty}(||\nabla\dot{S}_{[\frac{3k}{4}]+3}u||^2_{L^2}+||\nabla\dot{S}_{[\frac{3k}{4}]+3}B||^2_{L^2}),\notag
\end{align}
\begin{align}\label{b7}
|J_7|=&|\int_{\mathbb{R}^3}\sum_{l=k-1}^{[\frac{3}{4} k]}\tilde{\dot{\Delta}}_{l}B^k\cdot\nabla\dot{S}_kB\cdot\dot{\Delta}_lu^kdx|\\
\leq&||\nabla\dot{S}_kB||_{L^\infty}\sum_{l=k-1}^{[\frac{3k}{4}]}||\tilde{\dot{\Delta}}_l B^k||_{L^2}||\dot{\Delta}_lu^k||_{L^2}\notag\\
\leq& C2^k||\dot{S}_kB||_{L^\infty}\sum_{l=k-1}^{[\frac{3k}{4}]}2^{-2l}||\nabla\dot{\Delta}_l u||_{L^2}||\nabla\tilde{\dot{\Delta}}_lB||_{L^2}\notag\\
\leq& C2^{-k}||\dot{S}_kB||_{L^\infty}\sum_{l=k-1}^{[\frac{3k}{4}]}2^{2(k-l)}||\nabla\dot{\Delta}_l u||_{L^2}||\nabla\tilde{\dot{\Delta}}_lB||_{L^2}\notag\\
\leq&C2^{-k}||\dot{S}_kB||_{L^\infty}(||\nabla\dot{S}_{[\frac{3k}{4}]+3}u||^2_{L^2}+||\nabla\dot{S}_{[\frac{3k}{4}]+3}B||^2_{L^2}),\notag
\end{align}
\begin{align}\label{b8}
|J_8|=&|\int_{\mathbb{R}^3}\sum_{l=k-1}^{[\frac{3}{4} k]}\tilde{\dot{\Delta}}_{l}B^k\cdot\nabla\dot{S}_ku\cdot\dot{\Delta}_lB^kdx|\\
\leq&||\nabla\dot{S}_ku||_{L^\infty}\sum_{l=k-1}^{[\frac{3k}{4}]}||\tilde{\dot{\Delta}}_l B^k||_{L^2}||\dot{\Delta}_lB^k||_{L^2}\notag\\
\leq& C2^k||\dot{S}_ku||_{L^\infty}\sum_{l=k-1}^{[\frac{3k}{4}]}2^{-2l}||\nabla\dot{\Delta}_l B||_{L^2}||\nabla\tilde{\dot{\Delta}}_lB||_{L^2}\notag\\
\leq& C2^{-k}||\dot{S}_ku||_{L^\infty}\sum_{l=k-1}^{[\frac{3k}{4}]}2^{2(k-l)}||\nabla\dot{\Delta}_l B||_{L^2}||\nabla\tilde{\dot{\Delta}}_lB||_{L^2}\notag\\
\leq&C2^{-k}||\dot{S}_ku||_{L^\infty}||\nabla\dot{S}_{[\frac{3k}{4}]+3}B||^2_{L^2}.\notag
\end{align} Substituting the estimates \eqref{b1}-\eqref{b8} into \eqref{+11}, we deduce that
\begin{align}\label{BDI}
&\int_{\mathbb{R}^3}(|\nabla u|^2+|\nabla B|^2)dx\\
\leq& C\liminf_{k\to-\infty}2^{-k}(||\dot{S}_{k}u||_{L^\infty}+||\dot{S}_{k}B||_{L^\infty})(||\nabla\dot{S}_{[\frac{3k}{4}]+3}B||^2_{L^2}+||\nabla\dot{S}_{[\frac{3k}{4}]+3}u||^2_{L^2})\notag
\end{align}where we have used the fact that $||\nabla\dot{S}_{k+3}B||^2_{L^2}\leq C||\nabla\dot{S}_{[\frac{3k}{4}]+3}B||^2_{L^2}$ as $k\to-\infty$.

 Notice that $(||\nabla u||^2_{L^2}+||\nabla B||^2_{L^2})=D(u,B)<\infty$, one can see that
\begin{align}\label{zero}
\lim_{k\to-\infty}(||\nabla\dot{S}_{[\frac{3k}{4}]+3}u||^2_{L^2}+||\nabla\dot{S}_{[\frac{3k}{4}]+3}B||^2_{L^2})=0.
\end{align}From \eqref{BDI} and \eqref{zero}, we deduce that if $\liminf_{k\to-\infty}2^{-k}(||\dot{S}_{k}u||_{L^\infty}+||\dot{S}_{k}B||_{L^\infty})<\infty$, then it follows
\begin{align*}
\int_{\mathbb{R}^3}(|\nabla u|^2+|\nabla B|^2)dx=0.
\end{align*}
From the fact $||f||_{L^6(\mathbb{R}^3)}\leq C||\nabla f||_{L^2(\mathbb{R}^3)}$, we thus conclude $u=B=0$ and obtain the first conclusion in Theorem \ref{main}.

Furthermore, noticing that
\begin{align*}
&2^{-k}||\dot{S_k}f||_{L^\infty}\leq2^{-k}\sum_{l\leq {k-1}}||\dot{\Delta}_lf||_{L^\infty}\\
\leq& 2^{-k}\sum_{l\leq {k-1}}2^l||\dot{S}_kf||_{\dot{B}^{-1}_{\infty,\infty}}\leq C||\dot{S}_kf||_{\dot{B}^{-1}_{\infty,\infty}},\notag
\end{align*}we have that
\begin{align}\label{besov}
\liminf_{k\to-\infty}2^{-k}(||\dot{S}_{k}u||_{L^\infty}+||\dot{S}_{k}B||_{L^\infty})\leq C\liminf_{k\to-\infty}(||\dot{S}_{k}u||_{\dot{B}^{-1}_{\infty,\infty}}+||\dot{S}_{k}B||_{\dot{B}^{-1}_{\infty,\infty}}).
\end{align} This means that if $\liminf_{k\to-\infty}(||\dot{S}_{k}u||_{\dot{B}^{-1}_{\infty,\infty}}+||\dot{S}_{k}B||_{\dot{B}^{-1}_{\infty,\infty}})<\infty$ then $u=B=0$. We
obtain the second conclusion in Theorem \ref{main} and
complete the proof of Theorem \ref{main}.
\qquad $\hfill\Box$

Next, we start the proof of Theorem \ref{main1}.

\emph{\bf{Proof of Theorem \ref{main1}.}} We will establish new estimates for the right hand side terms in \eqref{+11} based on inequalities \eqref{Ber1}-\eqref{+10} and the embedding theorem $||f||_{L^6(\mathbb{R}^3)}
\leq C||\nabla f||_{L^2(\mathbb{R}^3)}$. We will apply similar computations of \eqref{b1}-\eqref{b8} to control these terms, for each term, we will replace the $L^\infty\times L^2\times L^2$ estimates by the $L^3\times L^6\times L^2$ estimates. The main differences come from the terms $J_6$ and $J_7$. For the interactions between the velocity field $u$ and the magnetic field $B$ at the high frequency region, we will apply the well-known fact $||(\sum_j|\dot{\Delta}_jf|^2)^\frac{1}{2}||_{L^p}=||f||_{L^p}$ with $1<p<\infty$ to establish the $L^3\times L^6\times L^2$ estimates for $J_6$ and $J_7$ as follows:
\begin{align}\label{nb6}
|J_6|=&|\int_{\mathbb{R}^3}\sum_{l=k-1}^{[\frac{3}{4} k]}\tilde{\dot{\Delta}}_{l}u^k\cdot\nabla\dot{S}_kB\cdot\dot{\Delta}_lB^kdx|\\
\leq&\int_{\mathbb{R}^3}|\nabla\dot{S}_kB|(\sum_{l=k-1}^{[\frac{3}{4} k]}|\tilde{\dot{\Delta}}_{l}u^k|^2)^\frac{1}{2}(\sum_{l=k-1}^{[\frac{3}{4} k]}|{\dot{\Delta}}_{l}B^k|^2)^\frac{1}{2}dx\notag\\
\leq&||\nabla\dot{S}_kB||_{L^2}||(\sum_{l=k-1}^{[\frac{3}{4} k]}|\tilde{\dot{\Delta}}_{l}u^k|^2)^\frac{1}{2}||_{L^3}||(\sum_{l=k-1}^{[\frac{3}{4} k]}|{\dot{\Delta}}_{l}B^k|^2)^\frac{1}{2}||_{L^6}\notag\\
\leq&||\nabla\dot{S}_kB||_{L^2}||\dot{S}_{[\frac{3k}{4}]+3}u||_{L^3}||\dot{S}_{[\frac{3k}{4}]+1}B||_{L^6}\notag\\
\leq& C||\dot{S}_{[\frac{3k}{4}]+3}u||_{L^3}||\nabla\dot{S}_kB||_{L^2}||\nabla\dot{S}_{[\frac{3k}{4}]+1}B||_{L^2}\notag\\
\leq& C||\dot{S}_{[\frac{3k}{4}]+3}u||_{L^3}(||\nabla\dot{S}_kB||^2_{L^2}+||\nabla\dot{S}_{[\frac{3k}{4}]+1}B||^2_{L^2}),\notag
\end{align}
\begin{align}\label{nb7}
|J_7|=&|\int_{\mathbb{R}^3}\sum_{l=k-1}^{[\frac{3}{4} k]}\tilde{\dot{\Delta}}_{l}B^k\cdot\nabla\dot{S}_kB\cdot\dot{\Delta}_lu^kdx|\\
\leq&\int_{\mathbb{R}^3}|\nabla\dot{S}_kB|(\sum_{l=k-1}^{[\frac{3}{4} k]}|\tilde{\dot{\Delta}}_{l}B^k|^2)^\frac{1}{2}(\sum_{l=k-1}^{[\frac{3}{4} k]}|{\dot{\Delta}}_{l}u^k|^2)^\frac{1}{2}dx\notag\\
\leq&||\nabla\dot{S}_kB||_{L^2}||(\sum_{l=k-1}^{[\frac{3}{4} k]}|\tilde{\dot{\Delta}}_{l}B^k|^2)^\frac{1}{2}||_{L^6}||(\sum_{l=k-1}^{[\frac{3}{4} k]}|{\dot{\Delta}}_{l}u^k|^2)^\frac{1}{2}||_{L^3}\notag\\
\leq&||\nabla\dot{S}_kB||_{L^2}||\dot{S}_{[\frac{3k}{4}]+3}B||_{L^6}||\dot{S}_{[\frac{3k}{4}]+1}u||_{L^3}\notag\\
\leq& C||\dot{S}_{[\frac{3k}{4}]+1}u||_{L^3}||\nabla\dot{S}_kB||_{L^2}||\nabla\dot{S}_{[\frac{3k}{4}]+3}B||_{L^2}\notag\\
\leq& C||\dot{S}_{[\frac{3k}{4}]+1}u||_{L^3}(||\nabla\dot{S}_kB||^2_{L^2}+||\nabla\dot{S}_{[\frac{3k}{4}]+3}B||^2_{L^2}),\notag
\end{align}For the remaining terms $J_1-J_5$ and $J_8$, we will use the similar computations in the proof of Theorem \ref{main} to establish the $L^3\times L^6\times L^2$ estimates.
\begin{align}\label{nb1}
&|J_1|\\=&|\int_{\mathbb{R}^3}
\sum_{l=k-1}^{k+1}\sum_{l'=l-2}^{k-1}\dot{S}_{l-2}\dot{S}_{k}u\cdot\nabla\dot{\Delta}_{l'}u\cdot\dot{\Delta}_lu^kdx+\int_{\mathbb{R}^3}\sum_{l=k-3}^{k}
\dot{\Delta}_l\dot{S}_ku\cdot\nabla\dot{S}_ku\cdot\tilde{\dot{\Delta}}_lu^kdx|\notag\\
\leq&\sum_{l=k-1}^{k+1}\sum_{l'=l-2}^{k-1}||\dot{S}_{l-2}\dot{S}_ku||_{L^3}||\nabla\dot{\Delta}_{l'}u||_{L^2}||\dot{\Delta}_lu^k||_{L^6}
+\sum_{l=k-3}^{k}||\nabla\dot{S}_{k}u||_{L^2}||\dot{\Delta}_{l}\dot{S}_ku||_{L^3}||\tilde{\dot{\Delta}}_lu^k||_{L^6}\notag\\
\leq&C||\dot{S}_ku||_{L^3}(\sum_{l=k-1}^{k+1}\sum_{l'=l-2}^{k-1}||\nabla\dot{\Delta}_{l'}u||_{L^2}||\nabla\dot{\Delta}_lu||_{L^2}
+\sum_{l=k-3}^{k}||\nabla\dot{S}_{k}u||_{L^2}||\nabla\tilde{\dot{\Delta}}_lu||_{L^2})\notag\\
\leq& C||\dot{S}_ku||_{L^3}(||\nabla\dot{S}_{k}u||^2_{L^2}+||\nabla\dot{S}_{k+3}u||^2_{L^2}),\notag
\end{align}
\begin{align}\label{nb2}
&|J_2|\\=&|\int_{\mathbb{R}^3}\sum_{l=k-1}^{k+1}\sum_{l'=l-2}^{k-1}\dot{S}_{l-2}\dot{S}_{k}u\cdot\nabla\dot{\Delta}_{l'}B\cdot\dot{\Delta}_lB^kdx
+\int_{\mathbb{R}^3}\sum_{l=k-3}^{k}\dot{\Delta}_l\dot{S}_ku\cdot\nabla\dot{S}_kB\cdot\tilde{\dot{\Delta}}_lB^kdx|\notag\\
\leq&\sum_{l=k-1}^{k+1}\sum_{l'=l-2}^{k-1}||\dot{S}_{l-2}\dot{S}_ku||_{L^3}||\nabla\dot{\Delta}_{l'}B||_{L^2}||\dot{\Delta}_lB^k||_{L^6}
+\sum_{l=k-3}^{k}||\nabla\dot{S}_{k}B||_{L^2}||\dot{\Delta}_{l}\dot{S}_ku||_{L^3}||\tilde{\dot{\Delta}}_lB^k||_{L^6}\notag\\
\leq&C||\dot{S}_ku||_{L^3}(\sum_{l=k-1}^{k+1}\sum_{l'=l-2}^{k-1}||\nabla\dot{\Delta}_{l'}B||_{L^2}||\nabla\dot{\Delta}_lB||_{L^2}
+\sum_{l=k-3}^{k}||\nabla\dot{S}_{k}B||_{L^2}||\nabla\tilde{\dot{\Delta}}_lB||_{L^2})\notag\\
\leq& C||\dot{S}_ku||_{L^3}(||\nabla\dot{S}_{k}B||^2_{L^2}+||\nabla\dot{S}_{k+3}B||^2_{L^2}),\notag
\end{align}
\begin{align}\label{nb3}
&|J_3|\\=&|\int_{\mathbb{R}^3}\sum_{l=k-1}^{k+1}\sum_{l'=l-2}^{k-1}\dot{S}_{l-2}\dot{S}_{k}B\cdot\nabla\dot{\Delta}_{l'}B\cdot\dot{\Delta}_lu^kdx
+\int_{\mathbb{R}^3}\sum_{l=k-3}^{k}\dot{\Delta}_l\dot{S}_kB\cdot\nabla\dot{S}_kB\cdot\tilde{\dot{\Delta}}_lu^kdx|\notag\\
\leq&\sum_{l=k-1}^{k+1}\sum_{l'=l-2}^{k-1}||\dot{S}_{l-2}\dot{S}_kB||_{L^6}||\nabla\dot{\Delta}_{l'}B||_{L^2}||\dot{\Delta}_lu^k||_{L^3}
+\sum_{l=k-3}^{k}||\nabla\dot{S}_{k}B||_{L^2}||\dot{\Delta}_{l}\dot{S}_kB||_{L^6}||\tilde{\dot{\Delta}}_lu^k||_{L^3}\notag\\
\leq&C||\dot{S}_{k+3}u||_{L^3}(\sum_{l=k-1}^{k+1}\sum_{l'=l-2}^{k-1}||\nabla\dot{S}_kB||_{L^2}||\nabla\dot{\Delta}_{l'}B||_{L^2}
+\sum_{l=k-3}^{k}||\nabla\dot{S}_{k}B||_{L^2}||\nabla\dot{\Delta}_{l}B||_{L^2})\notag\\
\leq& C||\dot{S}_{k+3}u||_{L^3}(||\nabla\dot{S}_{k}B||^2_{L^2}+||\nabla\dot{S}_{k+3}B||^2_{L^2}),\notag
\end{align}
\begin{align}\label{nb4}
&|J_4|\\=&|\int_{\mathbb{R}^3}\sum_{l=k-1}^{k+1}\sum_{l'=l-2}^{k-1}\dot{S}_{l-2}\dot{S}_{k}B\cdot\nabla\dot{\Delta}_{l'}u\cdot\dot{\Delta}_lB^kdx
+\int_{\mathbb{R}^3}\sum_{l=k-3}^{k}\dot{\Delta}_l\dot{S}_kB\cdot\nabla\dot{S}_ku\cdot\tilde{\dot{\Delta}}_lB^kdx|\notag\\
\leq&\sum_{l=k-1}^{k+1}\sum_{l'=l-2}^{k-1}||\dot{S}_{l-2}\dot{S}_kB||_{L^6}||\nabla\dot{\Delta}_{l'}u||_{L^3}||\dot{\Delta}_lB^k||_{L^2}
+\sum_{l=k-3}^{k}||\nabla\dot{S}_{k}u||_{L^3}||\dot{\Delta}_{l}\dot{S}_kB||_{L^6}||\tilde{\dot{\Delta}}_lB^k||_{L^2}\notag\\
\leq&C||\dot{S}_{k}u||_{L^3}(\sum_{l=k-1}^{k+1}\sum_{l'=l-2}^{k-1}||\nabla\dot{S}_kB||_{L^2}2^{l'}||\dot{\Delta}_lB||_{L^2}
+\sum_{l=k-3}^{k}||\nabla\dot{\Delta}_{l}B||_{L^2}2^k||\tilde{\dot{\Delta}}_lB||_{L^2})\notag\\
\leq&C||\dot{S}_{k}u||_{L^3}(\sum_{l=k-1}^{k+1}\sum_{l'=l-2}^{k-1}||\nabla\dot{S}_kB||_{L^2}||\nabla\dot{\Delta}_lB||_{L^2}
+\sum_{l=k-3}^{k}||\nabla\dot{\Delta}_{l}B||_{L^2}||\nabla\tilde{\dot{\Delta}}_lB||_{L^2})\notag\\
\leq& C||\dot{S}_{k}u||_{L^3}(||\nabla\dot{S}_{k}B||^2_{L^2}+||\nabla\dot{S}_{k+3}B||^2_{L^2}),\notag
\end{align}
\begin{align}\label{nb5}
|J_5|=&|\int_{\mathbb{R}^3}\sum_{l=k-1}^{[\frac{3k}{4}]}\dot{\Delta}_l u^k\cdot\nabla \dot{S}_ku\cdot\tilde{\dot{\Delta}}_lu^kdx|\\
\leq&C||\dot{S}_ku||_{L^3}\sum_{l=k-1}^{[\frac{3k}{4}]}2^k||\dot{\Delta}_l u||_{L^2}||\tilde{\dot{\Delta}}_lu||_{L^6}\notag\\
\leq& C||\dot{S}_ku||_{L^3}\sum_{l=k-1}^{[\frac{3k}{4}]}||\nabla\dot{\Delta}_l u||_{L^2}||\nabla\tilde{\dot{\Delta}}_lu||_{L^2}\notag\\
\leq& C||\dot{S}_ku||_{L^3}||\nabla\dot{S}_{[\frac{3k}{4}]+3}u||^2_{L^2},\notag
\end{align}
\begin{align}\label{nb8}
|J_8|=&|\int_{\mathbb{R}^3}\sum_{l=k-1}^{[\frac{3}{4} k]}\tilde{\dot{\Delta}}_{l}B^k\cdot\nabla\dot{S}_ku\cdot\dot{\Delta}_lB^kdx|\\
\leq&C||\dot{S}_ku||_{L^3}\sum_{l=k-1}^{[\frac{3k}{4}]}2^k||\tilde{\dot{\Delta}}_l B^k||_{L^6}||\dot{\Delta}_lB^k||_{L^2}\notag\\
\leq& C||\dot{S}_ku||_{L^3}\sum_{l=k-1}^{[\frac{3k}{4}]}2^{k-l}||\nabla\dot{\Delta}_l B||_{L^2}||\nabla\tilde{\dot{\Delta}}_lB||_{L^2}\notag\\
\leq& C||\dot{S}_ku||_{L^3}\sum_{l=k-1}^{[\frac{3k}{4}]}||\nabla\dot{\Delta}_l B||_{L^2}||\nabla\tilde{\dot{\Delta}}_lB||_{L^2}\notag\\
\leq& C||\dot{S}_{k}u||_{L^3}(||\nabla\dot{S}_{[\frac{3k}{4}]+1}B||^2_{L^2}+||\nabla\dot{S}_{[\frac{3k}{4}]+3}B||^2_{L^2}).\notag
\end{align}
Substituting \eqref{nb1}-\eqref{nb8} into \eqref{+11}, we obtain that
\begin{align}\label{nBDI}
&\int_{\mathbb{R}^3}(|\nabla u|^2+|\nabla B|^2)dx\\
\leq& C\liminf_{k\to-\infty}(||\dot{S}_ku||_{L^3}+||\dot{S}_{[\frac{3k}{4}]+1}u||_{L^3}+||\dot{S}_{[\frac{3k}{4}]+3}u||_{L^3})\times\notag\\
&(||\nabla\dot{S}_{k}u||^2_{L^2}+||\nabla\dot{S}_{k+3}u||^2_{L^2}+||\nabla\dot{S}_{k}B||^2_{L^2}+||\nabla\dot{S}_{k+3}B||^2_{L^2}\notag\\
&+||\nabla\dot{S}_{[\frac{3k}{4}]+1}u||^2_{L^2}+||\nabla\dot{S}_{[\frac{3k}{4}]+1}B||^2_{L^2}+||\nabla\dot{S}_{[\frac{3k}{4}]+3}B||^2_{L^2}).\notag
\end{align}From \eqref{DI}, it is not difficult to see that
\begin{align}\label{vanish}
\lim_{k\to-\infty}&(||\nabla\dot{S}_{k}u||^2_{L^2}+||\nabla\dot{S}_{k+3}u||^2_{L^2}+||\nabla\dot{S}_{k}B||^2_{L^2}+||\nabla\dot{S}_{k+3}B||^2_{L^2}\\
&+||\nabla\dot{S}_{[\frac{3k}{4}]+1}u||^2_{L^2}+||\nabla\dot{S}_{[\frac{3k}{4}]+1}B||^2_{L^2}+||\nabla\dot{S}_{[\frac{3k}{4}]+3}B||^2_{L^2})\notag\\
=&0.\notag
\end{align}Collecting \eqref{nBDI} and \eqref{vanish}, we conclude that if $\liminf_{k\to-\infty}||\dot{S}_ku||_{L^3}<\infty$ then it follows that
\begin{align*}
\int_{\mathbb{R}^3}(|\nabla u|^2+|\nabla B|^2)dx=0.
\end{align*} Consequently, we deduce from the fact $||f||_{L^6(\mathbb{R}^3)}\leq C||\nabla f||_{L^2(\mathbb{R}^3)}$ that $u=B=0$. This completes the proof of Theorem \ref{main1}.
\qquad $\hfill\Box$

\section*{Acknowledgments}
The author is supported by the Construct Program of the Key Discipline in Hunan Province and NSFC (Grant No. 11871209), and the Hunan Provincial NSF (No. 2022JJ10033)
\\

\end{document}